\documentclass[11pt]{article}

\usepackage[margin=1in]{geometry}
\usepackage{amsmath,amssymb,amsthm}
\usepackage{mathtools}
\usepackage{graphicx}
\usepackage{hyperref}
\usepackage{cleveref}
\usepackage{cite}

\newtheorem{theorem}{Theorem}[section]
\newtheorem{lemma}[theorem]{Lemma}

\newtheorem{corollary}[theorem]{Corollary}

\theoremstyle{definition}
\newtheorem{definition}[theorem]{Definition}

\theoremstyle{remark}

\title{Clique counts under degree and matching constraints}

\author{
    Rajat Adak\\
    Department of Computer Science and Automation\\
    Indian Institute of Science, Bengaluru, India\\
    \texttt{rajatadak@iisc.ac.in}
}

\date{}

\begin{document}

\maketitle

\begin{abstract}
We study the maximum number of $t$-cliques in graphs with
maximum degree at most $D$ and matching number at most $M$,
with no restriction on the number of vertices. For $t\ge 3$,
we determine this maximum whenever $D\ge 2M-1$ and
characterize all extremal graphs when the maximum is
positive. Up to isolated vertices, these are complete or
complete split graphs for $D\ge 2M$, while $K_{2M}$ is
uniquely extremal when $D=2M-1$. We also derive a general
upper bound from the Chv\'atal-Hanson edge theorem and
refine it for even $D$. For fixed even $D$ and
$3\le t\le D+1$, every extremal graph consists, apart from
isolated vertices, of copies of $K_{D+1}$ together with a
remainder whose order is bounded independently of $M$.
When $D/2$ divides $M$, this remainder is empty and the
extremal graph is unique. This structure yields an eventual
recurrence, reducing the determination of all extremal
values and graphs to finitely many matching bounds.
\end{abstract}


\section{Introduction}
Extremal graph theory studies how large a specified subgraph count can be under certain structural restrictions. For graphs of fixed order, Wang~\cite{wang2020shifting} extended the Erd\H{o}s-Gallai matching theorem to clique counts, while Chao and Dong~\cite{Chao_Dong_2022} gave a unified proof of the Gan-Loh-Sudakov conjecture for clique counts under a maximum-degree constraint. Under the simultaneous restrictions of bounded maximum degree and matching number, Chv\'atal and Hanson~\cite{chvatal1976degrees} determined the maximum number of $K_2$'s, that is, edges. Balachandran and Khare~\cite{balachandran2009graphs} subsequently gave a structural proof of the same formula. In this paper, we study the corresponding problem for larger cliques under both constraints, without imposing any restriction on the order of the graph.

All graphs considered are finite and simple. Let $\Delta(G)$ and $\nu(G)$ denote the maximum degree and matching number of a graph $G$, respectively, and let $k_t(G)$ be its number of $t$-vertex cliques, counted by their vertex sets. We consider
\[F_t(D,M)=\max\{k_t(G):\Delta(G)\le D,\ \nu(G)\le M\}.\]
The endpoints of a maximum matching form a vertex cover of size at most $2M$, so $|E(G)|\le 2DM$. After isolated vertices are deleted, this also gives $|V(G)|\le 4DM$, and hence the maximum exists. Clearly, $F_t(D,0)=0$ for $t\ge 2$. We use the convention that $\binom{n}{j}=0$ for $n\ge 0$ whenever $j<0$ or $j>n$.

The parameters $\Delta(G)<d$ and $\nu(G)<m$ used by Balachandran and Khare~\cite{balachandran2009graphs} correspond to $D=d-1$ and $M=m-1$. For larger cliques, the edge maximum alone need not identify the extremal graphs.

In \Cref{thm1,thm2}, we determine $F_t(D,M)$ for $D\ge 2M-1$ and characterize the extremal graphs. We also give a general upper bound and improve it for even $D$
(\Cref{thm3,thm4}). In the even-degree case, \Cref{thm5} gives a bounded-remainder description of the extremal graphs, from which \Cref{cor_eventual_structure} derives an eventual recurrence.

Let $S_{D,M}=K_M\vee \overline{K_{D+1-M}},$
where $\vee$ denotes the join. If $D\ge 2M$, both $K_{2M+1}$ and $S_{D,M}$ have matching number $M$ and maximum degree at most $D$. Put
\begin{equation}\label{clique_count}
A_t(M)=\binom{2M+1}{t},\qquad
B_t(D,M)=\binom{M}{t}+(D+1-M)\binom{M}{t-1}.
\end{equation}
Clearly, $A_t(M) = k_t(K_{2M+1})$ and $B_t(D,M) = k_t(S_{D,M})$.

Note that the cases when $t > 2M+1$ or $t > D+1$ are not interesting. This is because a clique of order at least $2M+2$ contains $M+1$ matching edges itself, and on the other hand $K_{D+2}$ contains vertices of degree $D+1$. Thus, in either case, $F_t(D,M) = 0$. 

\section{Tight bound when \texorpdfstring{$D\ge 2M$}{D >= 2M}}

As our first result we determine $F_t(D,M)$ throughout the range $D\ge 2M$, in which the largest clique allowed by the matching constraint is also allowed by the degree constraint. 
\begin{theorem}\label{thm1}
Let $M\ge 1$, $D\ge 2M$, and $3\le t\le 2M+1$ be integers. Then
\begin{equation}\label{eqn_main}
F_t(D,M)=\max\{A_t(M),B_t(D,M)\}.
\end{equation}
The only extremal graph is $K_{2M+1}$ if $A_t(M)>B_t(D,M)$, and $S_{D,M}$ if $B_t(D,M)>A_t(M)$, up to isolated vertices. If the two values are equal, these two graphs are the only extremal graphs.
\end{theorem}

Wang~\cite{wang2020shifting} proved that the maximum number of $K_t$'s in an $n$-vertex graph with matching number at most $M$ is
\begin{equation*}
\max\left\{\binom{2M+1}{t},\binom{M}{t}+(n-M)\binom{M}{t-1}\right\}\qquad(n\ge 2M+1).
\end{equation*}
Our formula has the same form with $n$ replaced by $D+1$. But this substitution is not a consequence of Wang's theorem since a graph of maximum degree $D$ may have arbitrarily many vertices.

We first establish two counting lemmas for the proof of \Cref{thm1}. 

\begin{definition} For a set $S\subseteq V(G)$, let $c_t(G,S,{\le j})$ and $c_t(G,S,{\ge j})$ denote the numbers of $t$-cliques having at most $j$ vertices outside $S$ and at least $j$ vertices outside $S$, respectively. \end{definition}
\begin{lemma}\label[lemma]{lem1}
Let $S\subseteq V(G)$ with $|S|=s$, and suppose that $\Delta(G)\le D$ and $D\ge 2s$. For every $t\ge 3$,
\begin{equation}\label{eqn_lem1}
c_t(G,S,{\le 1})\le \binom{s}{t}+(D+1-s)\binom{s}{t-1}.
\end{equation}
If $s\ge t-1$, equality implies that $G[S]=K_s$, every vertex of $S$ has degree $D$, and every vertex outside $S$ that has a neighbour in $S$ is adjacent to all of $S$.
\end{lemma}

\begin{proof}
If $s<t-1$, both sides of \eqref{eqn_lem1} are zero. Thus we will assume $s\ge t-1$, and for $v\in S$ put $a_v=|N_G(v)\cap S|$. Let $c_0$ and $c_1$ be the numbers of $t$-cliques containing zero and one vertex outside $S$, respectively. Counting the cliques containing each vertex of $S$ gives
\[tc_0\le \sum_{v\in S}\binom{a_v}{t-1},\qquad (t-1)c_1\le \sum_{v\in S}(D-a_v)\binom{a_v}{t-2}.\]
For the second inequality, fix $v\in S$ and a neighbour $u\notin S$. The other $t-2$ vertices must lie in $N_G(v)\cap S$, giving at most $\binom{a_v}{t-2}$ choices. There are at most $D-a_v$ choices for $u$. 
Define 
\begin{equation*}
    f(a)=\frac1t\binom{a}{t-1}+\frac{D-a}{t-1}\binom{a}{t-2}.
\end{equation*}
Therefore, we get 
$c_t(G,S, \le 1)\le \sum_{v\in S} f(a_v).$
For $0\le a\le s-2$, Pascal's identity gives 
\[f(a+1)-f(a)=\frac{D-a-1}{t-1}\binom{a}{t-3}-\frac{1}{t(t-1)}\binom{a}{t-2}.\]
Clearly this difference is zero if $a<t-3$. Otherwise the difference is
\[\frac{1}{t-1}\binom{a}{t-3}\left(D-a-1-\frac{a-t+3}{t(t-2)}\right).\]
 Recall that $D \geq 2s, a \leq s-2$ and $t \geq 3$. Therefore, the bracket above evaluates to at least $s+1-(s-2)/3>0$, and therefore the difference is positive. Thus $f(a)\le f(s-1)$, with equality only at $a=s-1$. Therefore,
\[\sum_{v\in S}f(a_v)\le sf(s-1)=\binom{s}{t}+(D+1-s)\binom{s}{t-1}.\]

If equality holds then every $a_v=s-1$, so $G[S]=K_s$. Since $\binom{s-1}{t-2}>0$, equality in the bound for $c_1$ also forces every $v\in S$ to have exactly $D-a_v = D-s+1$ neighbours outside $S$. Since $G[S]$ is a clique of order $s$, every vertex in $S$ has degree $(s-1) + (D-s+1) = D$. For each outside neighbour $u$ of $v$, all $\binom{a_v}{t-2} = \binom{s-1}{t-2}$ possible choices must occur. As $1\le t-2\le s-1$, this requires $u$ to be adjacent to every vertex of $S\setminus\{v\}$. Therefore the equality conditions follow.
\end{proof}

\begin{lemma}\label[lemma]{lem_2}
Let $r\in \mathbb{Z}_{\geq 0}$, and let $H$ be a graph with $H_1,\ldots,H_q$ as its components. Suppose that
$\sum_{i=1}^q\left\lfloor\frac{|V(H_i)|}{2}\right\rfloor\le r.$
Then, for every $j\ge 2$,
\begin{equation}\label{eqn_lem2}
k_j(H)\le \binom{2r+1}{j}.
\end{equation}
If $2\le j\le 2r+1$, equality holds if and only if $H=K_{2r+1}$, up to isolated vertices.
\end{lemma}

\begin{proof}
We can ignore isolated vertices in $H$. Put $r_i=\lfloor |V(H_i)|/2\rfloor$ for each remaining component, and let $R=\sum_i r_i\le r$. Each component has at most $2r_i+1$ vertices. Now take complete graphs of orders $2r_i+1$ and identify one vertex of each with a single common vertex $z$, keeping all other vertices distinct. The resulting graph $J$ has $2R+1$ vertices, and every clique of order at least two lies within one of these complete graphs. Therefore
\[k_j(H)\le \sum_i\binom{2r_i+1}{j}=k_j(J)\le \binom{2R+1}{j}\le \binom{2r+1}{j}.\]

Suppose that $2\le j\le 2r+1$ and equality holds in \eqref{eqn_lem2} and thus in the above chain as well. Since the right-hand side is positive we get $R=r$. Equality $k_j(J)=\binom{2r+1}{j}$ forces $J$ to be complete. If there are at least two nontrivial components in $H$, then we will get two vertices in $J$ each corresponding to one of these components, which will be non-adjacent, resulting in a contradiction. Hence there is only one nontrivial component of $H$. That component must have $2r+1$ vertices, and all its $j$-sets must be cliques, so it is $K_{2r+1}$. The converse is easy to check.
\end{proof}

 Let $o(X)$ denote the number of odd components in $X$. The Tutte-Berge formula states that
\[\nu(G)=\frac{1}{2}\min_{S\subseteq V(G)}
\left(|V(G)|+|S|-o(G-S)\right).\]
Let $\mathcal{C}(X)$ denote the set of all the components of $X$. The Tutte-Berge formula is equivalent to
\begin{equation}\label{eqn_TB}
\nu(G)=\min_{S\subseteq V(G)}\left(|S|+\sum_{C\in\mathcal C(G-S)}\left\lfloor\frac{|V(C)|}{2}\right\rfloor\right).
\end{equation}
Let $C_1,C_2,\dots,C_k$ be the odd components of $G-S$ and $C_{k+1},\dots,C_{l}$ be the even ones. Let $\left\lfloor\frac{|V(C_i)|}{2}\right\rfloor = r_i$ for all $i \in [l]$. Since $|V(G)|-|S|=2\sum_{i=1}^{l}r_i+k$ and $o(G-S)=k$, we have
\begin{align*}
\frac{|V(G)|+|S|-o(G-S)}{2}
&=|S|+\frac{|V(G)|-|S|-o(G-S)}{2}\\
&=|S|+\frac{2\sum_{i=1}^{l}r_i+k-k}{2}=|S|+\sum_{i=1}^{l}r_i.
\end{align*}
\begin{proof}[Proof of \Cref{thm1}] Let $G$ satisfy $\Delta(G)\le D$, $\nu(G)\le M$ and $D \geq 2M$, and choose $S$ to be a Tutte-Berge set for $G$, that is, a set attaining the minimum in \eqref{eqn_TB}. Put $s=|S|$. Then from \eqref{eqn_TB},
\begin{equation}\label{eqn1_thm1}
0\le s\le M,\qquad \sum_{C\in\mathcal C(G-S)}\left\lfloor\frac{|V(C)|}{2}\right\rfloor\le M-s.
\end{equation}
Fix $3\le t\le 2M+1$. The cliques with at most one vertex outside $S$ are bounded by \Cref{lem1}. From \eqref{eqn1_thm1}, $D\ge 2M\ge 2s$. If a clique uses $j\ge 2$ vertices outside $S$, those vertices lie in a single component of $G-S$ and those vertices form a $j$-clique. There are at most $\binom{s}{t-j}$ choices for its vertices in $S$. 
Using \Cref{lem1,lem_2}, \eqref{eqn1_thm1}, and Vandermonde's identity, we get
\begin{align}\label{eqn2_thm1}
k_t(G) = c_t(G,S,{\le 1}) + c_t(G,S,{\ge 2}) &\le \binom{s}{t}+(D+1-s)\binom{s}{t-1}
 +\sum_{j=2}^t\binom{s}{t-j}\binom{2(M-s)+1}{j}\nonumber\\
&=\binom{2M+1-s}{t}+(D-2M+s)\binom{s}{t-1}.
\end{align}
Put
\[
H_t(s)=\binom{2M+1-s}{t}+(D-2M+s)\binom{s}{t-1}.
\]
We now aim to maximize $H_t(s)$ for integers $0\le s\le M$. For $0\le s\le M-2$, Pascal's identity gives
\begin{align}
\Delta^2H_t(s)&=H_t(s+2)-2H_t(s+1)+H_t(s)\nonumber\\
&=\binom{2M-1-s}{t-2}
 +(D-2M+s+2)\binom{s}{t-3}
 +2\binom{s}{t-2}\ge 0.
\end{align}
Since $D\ge 2M$, $\Delta^2H_t(s)$ is nonnegative and thus $H_t$ is a convex sequence. Its endpoint values are $H_t(0)=A_t(M)$ and $H_t(M)=B_t(D,M)$, so
\begin{equation}\label{eqn3_thm1}
k_t(G)\le H_t(s)\le \left(1-\frac{s}{M}\right)A_t(M)+\frac{s}{M}B_t(D,M)\le \max\{A_t(M),B_t(D,M)\}.
\end{equation}
For $M\ge 2$, the second forward difference at $s=0$ satisfies
\[
\Delta^2H_t(0)\ge\binom{2M-1}{t-2}>0.
\] Thus, the forward differences are nondecreasing and not all equal. Consequently, the average of any proper initial segment of the forward differences is strictly smaller than their overall average. This implies that $H_t(s)$ lies strictly below the line segment joining its endpoint values $H_t(0)$ and $H_t(M)$ at every integer $0<s<M$. Hence, the middle inequality in \eqref{eqn3_thm1} is strict whenever $0<s<M$. For $M=1$, there are no interior integers, so the assertion is vacuously true.

Both endpoint constructions are feasible and attain the corresponding counts in \eqref{clique_count}, proving \eqref{eqn_main}. Now suppose $G$ attains the maximum. \eqref{eqn3_thm1} and the analysis above imply that $s=0$ or $s=M$.

If $s=0$, the bound on the component orders in \eqref{eqn1_thm1} and equality in \Cref{lem_2} imply that $G=K_{2M+1}$ apart from isolated vertices. Clearly, this case is possible only when $A_t(M)\ge B_t(D,M)$.

If $s=M$, from \eqref{eqn1_thm1}, every component of $G-S$ has order one, so $S$ is a vertex cover. The equality value $B_t(D,M)$ must be positive, which implies $t-1\le M$. Thus, we can apply the equality condition from \Cref{lem1}. Equality in \Cref{lem1} forces $G[S]=K_M$ and every vertex of $S$ to have degree $D$. Moreover, all vertices outside $S$ with a neighbour in $S$ have the same neighbourhood $S$. Therefore, there are exactly $D+1-M$ of these vertices and all the other vertices outside $S$ are isolated. Thus $G=S_{D,M}$ apart from isolated vertices. This case is possible only when $B_t(D,M)\ge A_t(M)$.\end{proof}
\section{The boundary case \texorpdfstring{$D=2M-1$}{D = 2M-1}}
In the previous section, we studied the clique bound $F_t(D,M)$ when $D\ge 2M$. Here we will study what happens for the boundary case of $D = 2M-1$. It is easy to check that $t>2M$ is not interesting, since a clique with order more than $2M$ contains vertices of degree at least $2M$. 
\begin{theorem}\label{thm2}
For $M\ge 2$ and $3\le t\le 2M$,
\begin{equation*}
F_t(2M-1,M)=\binom{2M}{t}.
\end{equation*}
$K_{2M}$ is the unique extremal graph, up to isolated vertices. For $t>2M$ the maximum is zero.
\end{theorem}
For a vertex $v$ in $G$, let $T_v$ denote the set of $t$-cliques meeting $N[v]$. Chao and Dong~\cite[Lemma~2]{Chao_Dong_2022} proved that for $t \geq 3$,
\begin{equation}\label{eqn_GaoDong}
\sum_{v\in V(G)}|T_v|\le \sum_{v\in V(G)}\binom{d_G(v)+1}{t}.
\end{equation}
We now establish a lemma as a consequence of the above bound.
\begin{lemma}\label[lemma]{lem3}
Let $M\ge 2$ and $3\le t\le 2M$. If $|V(G)|\le 2M+1$ and $\Delta(G)\le 2M-1$, then $k_t(G)\le\binom{2M}{t}$. Equality holds only for $K_{2M}$, with possibly one isolated vertex.
\end{lemma}

\begin{proof}
It is easy to check that the result holds when $|V(G)|\le 2M$, or when $G$ has an isolated vertex. Suppose that $|V(G)|=2M+1$ and that $G$ has no isolated vertex. Suppose for a contradiction that $k_t(G)\ge\binom{2M}{t}$. For a vertex $v$ of degree $d$, at most $\binom{2M-d}{t}$ cliques do not intersect with $N[v]$, so
\[|T_v|\geq \binom{2M}{t}-\binom{2M-d}{t}.\]
If $1\le d\le 2M-2$, convexity of the binomial coefficients gives
\begin{equation*}
    \binom{d+1}{t}+\binom{2M-d}{t} \le \binom{2M-1}{t}+\binom{2}{t} <\binom{2M}{t}.
\end{equation*}
Therefore, $|T_v|>\binom{d+1}{t}$ whenever $1\le d\le 2M-2$. If $d=2M-1$, we instead get $|T_v|\ge \binom{2M}{t}-\binom{1}{t} =\binom{d+1}{t}.$
Since we assumed $G$ has no isolated vertex, every vertex has degree at least one. Moreover, by the handshaking lemma, not all $2M+1$ vertices can have degree $2M-1$. Thus, at least one vertex has degree at most $2M-2$. Therefore, summing these inequalities over $V(G)$ contradicts \eqref{eqn_GaoDong}. 
Hence, $k_t(G)<\binom{2M}{t}$ when $|V(G)|=2M+1$ and $G$ has no isolated vertex.
\end{proof}
We now use \Cref{lem3} to obtain the desired bound and extremal graph for $D=2M-1$.
\begin{proof}[Proof of \Cref{thm2}]
Let $D=2M-1$ and let $G$ be feasible. We will look into two cases.

\noindent\textbf{Case 1.} There exists a vertex $v\in V(G)$ such that $\nu(G-v)<\nu(G)$.

\noindent Since $\nu(G)\le M$, we have $\nu(G-v)\le M-1$. Also, $\Delta(G-v)\le 2M-1$. We first show that
\begin{equation}\label{eqn_case1}
k_t(G-v)\le \binom{2M-1}{t}.
\end{equation}
Suppose first that $3\le t\le 2M-1$. We may apply \Cref{thm1} to $G-v$ with matching bound $M-1$ and degree bound $D=2M-1$. Writing $r=M-1$, we obtain
\begin{equation*}
k_t(G-v)\le
\max\left\{
\binom{2r+1}{t},
\binom{r}{t}+(r+2)\binom{r}{t-1}
\right\}.
\end{equation*}
We claim that the first term is at least the second. For $3\le t\le r+1$, Vandermonde's identity gives
\begin{equation*}
\binom{2r+1}{t}=\sum_{i=0}^{t}\binom{r+1}{i}\binom{r}{t-i}\ge \binom{r}{t}+(r+1)\binom{r}{t-1}
+\binom{r+1}{2}\binom{r}{t-2}>\binom{r}{t}+(r+2)\binom{r}{t-1},
\end{equation*}
where the last strict inequality follows from
$\binom{r+1}{2}\binom{r}{t-2}>\binom{r}{t-1}.$

For $t>r+1$, both terms in the second expression become zero, so the claimed comparison follows. This proves \eqref{eqn_case1} for $3\le t\le 2M-1$, and it is easy to check for $t \ge 2M$.

Since $D=2M-1$, we have $|N(v)| \leq 2M-1$. The $t$-cliques in $G$ can be partitioned based on whether they contain $v$. The number of $t$-cliques containing $v$ is therefore at most $\binom{2M-1}{t-1}$. By \eqref{eqn_case1}, the number of $t$-cliques not containing $v$ is at most $\binom{2M-1}{t}$. Thus,
\begin{equation*}
k_t(G) \le \binom{2M-1}{t}+\binom{2M-1}{t-1} =\binom{2M}{t}.
\end{equation*}
Suppose $k_t(G)=\binom{2M}{t}$. Therefore,
$k_{t-1}(G[N(v)])=\binom{2M-1}{t-1}.$
Since $2\le t-1\le 2M-1$ and $d(v)\le 2M-1$, this forces $d(v)=2M-1$ and $G[N(v)]=K_{2M-1}$. Thus, $G[N[v]] = K_{2M}$.

Every vertex of this $K_{2M}$ already has degree $2M-1$, so the degree bound $\Delta(G)\le 2M-1$ implies that this clique is a connected component of $G$. Moreover, it contains a matching of size $M$, thus it uses up the entire matching budget. Hence no other component can contain an edge. Thus all vertices outside this component are isolated.
\vspace{2mm}

\noindent\textbf{Case 2.} $\nu(G-v)=\nu(G)$ for every vertex $v \in V(G)$. 

\noindent Gallai's lemma, in the form used by Balachandran and Khare~\cite[Theorem~3]{balachandran2009graphs}, says that every nontrivial component $C_i$ is factor-critical and has $2r_i+1$ vertices, where $r_i=\nu(C_i)\ge 1$. Let $R=\sum_i r_i\le M$. If there is at most one nontrivial component, then, after removing isolated vertices, $G$ has at most $2M+1$ vertices. Therefore \Cref{lem3} gives $k_t(G)\le\binom{2M}{t}$. Equality would require $G=K_{2M}$, possibly with one isolated vertex, which is impossible since every nontrivial factor-critical component has odd order.

Now suppose that there are at least two nontrivial components. Let $C$ be one such component and $\nu(C) = x$, where $1\le x\le R-1$. By applying \Cref{lem_2} to $G\setminus C$ and using convexity, we get
\begin{align*}
k_t(G) = k_t(G\setminus C) + k_t(C) &\le \binom{2(R-x)+1}{t}+ \binom{2x+1}{t}\\
&\le \binom{2R-1}{t}+\binom{3}{t}\le \binom{2M-1}{t}+\binom{3}{t}
<\binom{2M}{t}.
\end{align*}
The final strict inequality follows from Pascal's identity and $\binom{2M-1}{t-1}>\binom{3}{t}$ for $3\le t\le 2M$. Since the inequality is strict in Case 2, the only extremal graph comes from Case 1, which is $K_{2M}$ with possible isolated vertices.

Conversely, $K_{2M}$ has maximum degree $2M-1$, matching number $M$, and exactly $\binom{2M}{t}$ copies of $K_t$. 
\end{proof}
\section{Generic upper bound and improvement when \texorpdfstring{$D$}{D} is even}

The previous sections determine $F_t(D,M)$ when $D\ge 2M-1$. We now give an upper bound valid for all degree and matching bounds, and then improve it when $D$ is even. The argument also gives information about the structure of the extremal graphs.

Chv\'atal and Hanson~\cite{chvatal1976degrees} proved that, for $D\ge 1$ and $M\ge 0$,
\begin{equation}\label{eqn_CH}
F_2(D,M)=DM+
\left\lfloor\frac D2\right\rfloor
\left\lfloor\frac{M}{\lceil D/2\rceil}\right\rfloor.
\end{equation}
We first use this result to bound the number of larger cliques.

\begin{theorem}\label{thm3}
Let $D\ge 1$, $M\ge 0$, and $2\le t\le D+1$ be integers. Every graph $G$ with $\Delta(G)\le D$ and $\nu(G)\le M$ satisfies
\begin{equation}\label{eqn_generic_bound}
k_t(G)\le U_t(D,M)=
\frac{\binom{D-1}{t-2}}{\binom{t}{2}}
\left(
DM+\left\lfloor\frac D2\right\rfloor
\left\lfloor\frac{M}{\lceil D/2\rceil}\right\rfloor
\right).
\end{equation}
Therefore, $F_t(D,M)\le\lfloor U_t(D,M)\rfloor$.
\end{theorem}

\begin{proof}
Fix an edge $uv\in E(G)$. Every $t$-clique containing $uv$ is obtained by choosing a $(t-2)$-clique in $N(u)\cap N(v)$. Since $|N(u)\cap N(v)|\le D-1$, the edge $uv$ lies in at most $\binom{D-1}{t-2}$ copies of $K_t$. Counting the pairs consisting of a $t$-clique and one of its edges gives
\[\binom{t}{2}k_t(G)\le
\binom{D-1}{t-2}|E(G)|\le
\binom{D-1}{t-2}F_2(D,M),\]
where the last inequality follows from \eqref{eqn_CH}. Dividing by $\binom{t}{2}$ proves \eqref{eqn_generic_bound}. Since $k_t(G)$ is an integer, the bound for $F_t(D,M)$ also follows.
\end{proof}
Putting $t=2$ in \eqref{eqn_generic_bound} gives exactly \eqref{eqn_CH}. For $t\ge 3$, the bound need not be sharp, since an edge need not lie in all $\binom{D-1}{t-2}$ possible copies of $K_t$. We will now keep track of some of this loss.

For $D\ge 2$ and $3\le t\le D+1$, put
\begin{equation}\label{eqn_count_constants}
a=\binom{D-1}{t-2},\qquad
b=\binom{t}{2},\qquad
c=\binom{D-2}{t-3}.
\end{equation}
Note that $a,b,c>0$ throughout this range.

\begin{lemma}\label[lemma]{lem4}
Let $\Delta(G)\le D$, where $D\ge 2$, and let $3\le t\le D+1$. Obtain $H$ from $G$ by deleting all components isomorphic to $K_{D+1}$ and all isolated vertices. Then
\begin{equation}\label{eqn_deficit}
bk_t(G)\le a|E(G)|-\frac c2|V(H)|.
\end{equation}
\end{lemma}

\begin{proof}
Call an edge $uv$ deficient if $|N(u)\cap N(v)|\le D-2$, and let $B$ be the set of all deficient edges. Every edge lies in at most $a$ copies of $K_t$, whereas a deficient edge lies in at most
\[\binom{D-2}{t-2}=a-c\]
such copies. Therefore, counting the number of $t$-cliques containing each edge gives
\[bk_t(G)\le a|E(G)|-c|B|.\]

We claim that every vertex of $H$ is incident with a deficient edge. Let $v\in V(H)$. If $d_G(v)<D$, then every edge incident with $v$ is deficient. Suppose instead that $d_G(v)=D$ and no edge incident with $v$ is deficient. For each $u\in N(v)$, this gives $|N(u)\cap N(v)|=D-1$, so $u$ is adjacent to every vertex of $N(v)\setminus\{u\}$. Thus $G[N[v]]=K_{D+1}$. Every vertex of this clique already has degree $D$, so it is a component of $G$. This contradicts the choice of $H$.

Since $H$ has no isolated vertices, the claim gives $2|B|\ge |V(H)|$. Substituting this into the preceding inequality proves \eqref{eqn_deficit}.
\end{proof}
Notice that \Cref{lem4} does not require $D$ to be even. Let $D=2h$, write $M=qh+r$, where $0\le r<h$, and put $A=\binom{D+1}{t}$. Define $f_0=0$ and
\begin{equation}\label{eqn_even_remainder}
f_r=\max\left\{
\binom{2r+1}{t},\,
\binom{r}{t}+(D+1-r)\binom{r}{t-1}
\right\}
\qquad(1\le r<h).
\end{equation}
These two expressions are the clique counts of the complete and complete split constructions introduced earlier, respectively. We write $qK_{D+1}$ for the disjoint union of $q$ copies of $K_{D+1}$.

\begin{theorem}\label{thm4}
Let $D=2h\ge 2$, $3\le t\le D+1$, and $M\ge 0$ be integers. Write $M=qh+r$, where $0\le r<h$. Then
\begin{equation}\label{eqn_even_bound}
qA+f_r\le F_t(D,M)\le
qA+\left\lfloor\frac{(Da-c)r}{b}\right\rfloor.
\end{equation}
\end{theorem}

\begin{proof}
For $r>0$, let $J_r$ be whichever of $K_{2r+1}$ and $S_{D,r}$ has more $t$-cliques, and let $J_0$ be the empty graph. Since $r<h$, we have $\Delta(J_r)\le D$ and $\nu(J_r)\le r$. Thus $qK_{D+1}\mathbin{\dot\cup}J_r$ is feasible and has $qA+f_r$ copies of $K_t$, proving the lower bound.

For the upper bound, let $G$ be feasible. We may assume that $\nu(G)=M$, otherwise, add $M-\nu(G)$ disjoint copies of $K_2$. Since $D\ge 2$ and $t\ge 3$, this does not change $k_t(G)$ or violate the degree bound.

Let $p$ be the number of components of $G$ isomorphic to $K_{D+1}$, and obtain $H$ by deleting these components and all isolated vertices. Each deleted $K_{D+1}$-component has matching number $h$, so $p\le q$ and
$\nu(H)=M-ph=(q-p)h+r\ge r.$
Therefore, $|V(H)|\ge 2r$. Also, \eqref{eqn_CH} gives
\begin{equation}\label{eqn_even_edges}
|E(G)|\le DM+h\left\lfloor\frac Mh\right\rfloor
=D(qh+r) + hq = qh(D+1)+Dr.
\end{equation}
Since $bA=ah(D+1)$, Lemma~\ref{lem4} gives
\begin{equation*}
bk_t(G) \le a|E(G)|-\frac c2|V(H)|\le a\bigl(qh(D+1)+Dr\bigr)-cr=bqA+(Da-c)r.
\end{equation*}
Dividing by $b$ and using integrality proves the upper bound.
\end{proof}
For the same parameters as above, \Cref{thm3} gives
\[F_t(D,M)\le qA+\left\lfloor\frac{Dar}{b}\right\rfloor.\]
Thus \eqref{eqn_even_bound} subtracts $cr/b$ from the remainder term before rounding. This is a strict improvement of the unrounded estimate when $r>0$, although the rounded bounds can coincide. When $r=0$, both bounds are attained by $qK_{D+1}$.

We now use the same argument to describe every extremal graph when $D$ is even. 

\begin{theorem}\label{thm5}
Let $D=2h\ge 2$, $3\le t\le D+1$, and $M\ge 0$ be integers, and put $r=M\bmod h$. Every graph attaining $F_t(D,M)$, after its isolated vertices are deleted, has the form
\begin{equation}\label{eqn_extremal_structure}
G=pK_{D+1}\mathbin{\dot\cup}H,
\end{equation}
where $H$ has no $K_{D+1}$-component and is extremal for $F_t(D,M-ph)$. Moreover,
\begin{equation}\label{eqn_remainder_bounds}
0\le M-ph\le \left\lfloor\frac{Dar-bf_r}{c}\right\rfloor,\qquad
|V(H)|\le\left\lfloor \frac{2(Dar-bf_r)}{c}\right\rfloor.
\end{equation}
Therefore, if $M=qh$, then
\begin{equation}\label{eqn_divisible}
F_t(D,M)=q\binom{D+1}{t},
\end{equation}
and $qK_{D+1}$ is the unique extremal graph, up to isolated vertices.
\end{theorem}

\begin{proof}
Write $M=qh+r$, where $0\le r<h$, and let $G$ be extremal. Delete all its $p$ components isomorphic to $K_{D+1}$ and all isolated vertices to obtain $H$. Since $\Delta(G)\le D$, every copy of $K_{D+1}$ must be a connected component of $G$, as each of its vertices already has $D$ neighbours within the clique and hence cannot have any neighbour outside it. Thus, $G=pK_{D+1}\mathbin{\dot\cup}H$, up to isolated vertices, and $H$ has no $K_{D+1}$-component. By \Cref{thm4}, $k_t(G)\ge qA+f_r$. Combining Lemma~\ref{lem4} with \eqref{eqn_even_edges} gives
\begin{equation}\label{eqn_order_bound}
\frac c2|V(H)|
\le a|E(G)|-bk_t(G)\le a\bigl(qh(D+1)+Dr\bigr)-b(qA+f_r)=Dar-bf_r.
\end{equation}
This proves $|V(H)|\le\lfloor \frac{2(Dar-bf_r)}{c}\rfloor$.

Since $\nu(H)\le M-ph$, the graph $H$ is feasible for $F_t(D,M-ph)$. If it were not extremal for this remaining matching bound, replacing it by a better feasible graph would increase $k_t(G)$ without violating either constraint. Thus $H$ is extremal for $F_t(D,M-ph)$.

It remains to bound $M-ph$. The inequality $ph\le\nu(G)\le M$ gives $M-ph\ge 0$. To obtain the upper bound even when $\nu(G)<M$, add $M-\nu(G)$ disjoint copies of $K_2$ to $G$, and denote the resulting graph by $G^+$. This graph is still feasible and extremal, has matching number $M$, and has the same $p$ components isomorphic to $K_{D+1}$. Its remainder $H^+$ therefore has matching number $M-ph$. Applying the order bound already proved in \eqref{eqn_order_bound} to $G^+$ gives
\[2(M-ph)\le |V(H^+)|\le \frac{2(Dar-bf_r)}{c}.\]
Hence $M-ph\le\lfloor\frac{Dar-bf_r}{c}\rfloor$, proving \eqref{eqn_remainder_bounds}.

Finally, if $r=0$, then $f_0=\frac{Dar-bf_r}{c}=0$. Thus $H$ is empty and $M-ph=0$, so $p=q$. The graph $qK_{D+1}$ attains \eqref{eqn_divisible}, proving the uniqueness assertion.
\end{proof}
\begin{corollary}\label[corollary]{cor_eventual_structure}
Let $D=2h\ge 2$, $3\le t\le D+1$, and $M\ge h$ be integers,
and put $r=M\bmod h$. If
$M>\left\lfloor\frac{Dar-bf_r}{c}\right\rfloor,$ then
\begin{equation}\label{eqn_eventual_recurrence}
F_t(D,M)=\binom{D+1}{t}+F_t(D,M-h).
\end{equation}
Moreover, up to isolated vertices, the extremal graphs for
$F_t(D,M)$ are exactly the graphs
\[
G=K_{D+1}\mathbin{\dot\cup}G',
\]
where $G'$ is extremal for $F_t(D,M-h)$.

Consequently, for fixed $D$ and $t$, determining $F_t(D,M)$
and its extremal graphs for all $M$ reduces to determining
them for finitely many values of $M$.
\end{corollary}

\begin{proof}
Let $G$ be extremal for $F_t(D,M)$, and let $p$ be its
number of components isomorphic to $K_{D+1}$. By
\Cref{thm5},
\[
M-ph\le\left\lfloor\frac{Dar-bf_r}{c}\right\rfloor<M.
\]
Therefore $p\ge 1$, so $G$ contains a $K_{D+1}$-component.
Delete one such component to obtain $G'$. Since
$\nu(K_{D+1})=h$ and matching number is additive over
components, we have
\[
\Delta(G')\le D,\qquad
\nu(G')=\nu(G)-h\le M-h.
\]
Thus $G'$ is feasible for $F_t(D,M-h)$, giving
\[
F_t(D,M)=k_t(G)
=\binom{D+1}{t}+k_t(G')
\le\binom{D+1}{t}+F_t(D,M-h).
\]

Conversely, adjoining a disjoint copy of $K_{D+1}$ to any
extremal graph for $F_t(D,M-h)$ gives a graph with maximum
degree at most $D$, matching number at most $M$, and
exactly $\binom{D+1}{t}+F_t(D,M-h)$ copies of $K_t$.
This proves \eqref{eqn_eventual_recurrence} and shows that
every graph constructed in this way is extremal.

For the graph $G'$ obtained by deleting a $K_{D+1}$-component
from an extremal graph $G$, the recurrence gives
\[
k_t(G')=F_t(D,M)-\binom{D+1}{t}=F_t(D,M-h).
\]
Hence $G'$ must itself be extremal, proving the claimed
characterization.

Finally, for fixed $D$ and $t$, the threshold
$\lfloor(Dar-bf_r)/c\rfloor$ depends only on the remainder
$r$. Subtracting $h$ preserves this remainder. Therefore,
we may apply the recurrence and the extremal characterization
repeatedly until the remaining matching bound is at most
this threshold or is less than $h$. Since $r$ has only
finitely many possible values, only finitely many matching
bounds remain to be considered.
\end{proof}

\bibliographystyle{plain}
\bibliography{references}

\end{document}